\documentclass[12pt,oneside,draft]{amsart}
\usepackage{amssymb,amsmath}

 

\newtheorem{theorem}[subsection]{Theorem}
\newtheorem*{theorem*}{Theorem}

\newtheorem{conjecture}[subsection]{Conjecture}
\theoremstyle{definition}
\newtheorem{definition}[subsection]{Definition}

\newtheorem{example}[subsection]{Example}
\theoremstyle{remark}
\newtheorem{remark}[subsection]{Remark}

\newcommand{\mt}[1]{\operatorname{#1}}

\newcommand{\QQ}{{\mathbb Q}}

\newcommand{\CC}{{\mathbb C}}

\newcommand{\PP}{{\mathbb P}}

\newcommand{\Diff}{\mt{Diff}}
\newcommand{\Exc}{\mt{Exc}}

\newcommand{\Supp}{\mt{Sup}}
\newcommand{\down}[1]{\llcorner #1 \lrcorner}
\newcommand{\up}[1]{\ulcorner #1 \urcorner}

\title{On exceptional terminal singularities}
\author{S.~A.~Kudryavtsev}

\date{}
\begin{document}
\begin{abstract} The first examples of exceptional terminal singularities
are constructed.
\end{abstract}
\maketitle

\section*{\bf {Introduction}}

The exceptional singularity study importance follows from the next observation.
\begin{enumerate}
\item If $(X\ni P)$ is an exceptional singularity then the linear system
$|-nK_X|$ is to have a "good" member for small $n$. Actually, we can take
$n \in \{1,2\}$ \cite[5.2]{Sh1} in two dimensional case and
$n \in \{1,2,3,4,6\}$ \cite[7.1]{Sh2} in three dimensional case.

\item Exceptional singularities are "bounded" and are to be classified.
\end{enumerate}

In this paper we construct the first examples of exceptional terminal
singularities.

\subsubsection*{Acknowledgements}
I am grateful to Professor V.A.Iskovskikh and Professor Yu.G.Prokhorov for
useful discussions and valuable remarks. The research was partially
supported by a grant 99-01-01132 from the Russian Foundation of Basic Research
and a grant INTAS-OPEN 97-2072.

\section{\bf Preliminaries}
All varieties are algebraic and are assumed to be defined over
$\CC$, complex number field. We will use the terminology and notation
of Log Minimal Model Program and the main properties of complements
\cite{Koetal}, \cite{PrLect}.

\begin{definition}
Let $(X/Z\ni P,D)$ be a contraction of varieties, where $D$ is a boundary.
Then a $\QQ$-{\it complement} of this contraction is an effective
$\QQ$-divisor $D'$ such that
$D'\ge D$, $K_X+D'$ is lc and $K_X+D'\sim_{\QQ} 0 $.
\end{definition}

\begin{definition}
Let $(X/Z\ni P,D)$ be a contraction of varieties, where $D$ is a boundary.
\begin{enumerate}
\item Assume that $Z$ is not a point (local case). Then $(X/Z\ni P,D)$
is said to be \textit{exceptional}
over $P$ if for any $\QQ$-complement of $K_X+D$ near the
fiber over $P$ there exists at most one (not necessarily
exceptional) divisor $E$ such that $a(E,D)=-1$.

\item Assume that $Z$ is a point (global case). Then $(X,D)$ is said to
be \textit{exceptional} if every $\QQ$-complement of $K_X+D$ is
klt.
\end{enumerate}
\end{definition}

\begin{example}\label{lc}
Let $(X \ni P)$ be a singularity. Suppose that there is an effective divisor
$H$ such that $(X,H)$ is lc and $\down{H} \ne 0$. Then the singularity is not
exceptional. Therefore three dimensional terminal singularity is not exceptional
because there is a divisor having only Du Val singularities in the anticanonical
linear system $|-K_X|$ \cite[6.4]{YPG}.
\end{example}

\begin{definition} \label{defplt}
Let $X$ be a normal lc  variety and let $f:Y \to X$ be a blow-up such that
the exceptional locus of $f$ contains only one irreducible divisor
$E\ (\Exc(f)=E)$. Then $f:(Y,E) \to X$ is called  {\it a purely log terminal
(plt) blow-up}, if $K_Y+E$ is plt and
$-E$ is $f$-ample.
\end{definition}

\begin{remark} \label{smain1}
If $(X \ni P)$ is klt singularity then there is a plt blow-up \cite[1.5]{Koud1}.
\end{remark}

\begin{theorem}\cite[4.9]{PrPlt}\label{CrExc}
Let $(X \ni P)$ be a klt singularity and let $f\colon (Y,E) \to X$ be a plt
blow-up of $P$. Then the following conditions are equivalent:
\begin{enumerate}
\item $(X \ni P)$ is an exceptional singularity;
\item $(E,\Diff_E(0))$ is an exceptional log variety.
\end{enumerate}
\end{theorem}

\section{\bf The examples of the exceptional terminal singularities}

\begin{theorem}
Let $(f=0,0)=(x^{a_1}_1+x^{a_2}_2+x^{a_3}_3+x^{a_4}_4+x^{a_5}_5=0,0)\\
\subset (\CC^5,0)$ be a four dimensional hypersurface singularities, where
$(a_1,\ldots,a_5)=$ $(2,3,11,17,19)$, $(2,3,11,17,23)$, $(2,3,11,17,25)$,
$(2,3,11,17,29)$, $(2,5,7,9,11)$, $(2,5,7,9,13)$.
Then they are terminal and exceptional.
\begin{proof}
Consider the first singularity. Let us prove that it is terminal. Since it is
given by non-degenerate polynomial then there exists embedded toric log
resolution \cite{Var}.
Therefore it is sufficient to prove that
$a_{\bf p}=\langle{\bf p},1\rangle - {\bf p}(f)-1 \ge 1$ for all
${\bf p}$, where
${\bf p}$ is non-zero vector with integral non-negative coordinates and
$\mbox{\bf p}(f) =\min_{x^m\in f}\langle\mbox{\bf p},m\rangle$.
The easy way to prove this statement is following.
Let $h(d)=\up{\frac d2}+\up{\frac d3}+\up{\frac d{11}}+\up{\frac
d{17}}+ \up{\frac d{19}}-d-1$. Let $d={\bf p}(f)$ then $h(d)\le a_{\bf p}$.
It is enough to check that
$h(d)\ge 1$ for all $1\le d\le 2\cdot 3\cdot11\cdot17\cdot19$.
The later is elementary to prove with the help of the computer program.
\par
The weighted blow-up of $\CC^5$ with weights proportional to
$(\frac12,\frac13,\frac1{11},\frac1{17},\frac1{19})$ induces a plt blow-up of
our singularity. The obtained log Fano variety
$(E,\Diff_E(0))$ is
\begin{equation*}
\begin{array}{cc}
\Big(\underbrace{\sum_{i=1}^5 x_i \subset \PP(1,1,1,1,1}),& \frac12H_1+\frac23H_2+\frac{10}{11}H_3+\frac{16}{17}H_4
+\frac{18}{19}H_5\Big),\\
\ \PP^3 &\\
\end{array}
\end{equation*}
where $H_i=\{x_i=0\}$.
Let $D=\sum d_iD_i+\sum h_iH_i$ be any $\QQ$-complement of this log variety. In
our case it is easy to check that
$d_i < 1$ and $h_i < 1$ for all $i$. If any component
$D_i$ of $D$ is a hyperplane section then we can easily show
$K_{\PP^3}+D$ to be klt.
Now let us prove that
$(\PP^3,D)$ is klt. This question is local. Therefore consider any chart
$\CC^3_{y_1,y_2,y_3}$.Let
$D=\Delta+D'$, where $\Delta$ is a sum of hyperplane sections and
$D'=\sum d_i\{f_i=0\}$. Let multiplicity of $f_i$ is equal to $n_i$.
By considering a change of coordinates of $\CC^3$ we can assume
without loss of generality that
$f_i=y^{n_i}_{1}+\cdots.$
Consider the deformation
$F'_t=\sum d_it^{-n_i}f_i(ty_1,t^2y_2,t^2y_3).$
For $t=0$ we have $D'_0=\{F'_0=0\}=\sum d_i\{y^{n_i}_1=0\}=\sum d_in_i\{y_1=0\}$
and for small
$t\ne 0$ we get $D'_t=\{F'_t=0\}=D'$.
By \cite[8.6]{Kollar1} it follows that
$c(\CC^3,\Delta + D'_0)\le c(\CC^3,\Delta + D'_t)=c(\CC^3,\Delta + D')$, where
$c$ is lc threshold.
Since the every component $\Supp(\Delta+D'_0)$ is a hyperplane section then
$(\CC^3,\Delta + D'_0)$ is klt. Therefore
$(\CC^3,\Delta + D')$ is also klt.
By criterion \ref{CrExc} the singularity is exceptional.
Divisor $\frac1{22}(\down{ 23\cdot\Diff_E(0)})$ gives
a 22-complement of minimal index.
The other singularities are considered similarly.
They have 24,34,34,22,28-complements of minimal index respectively.
\end{proof}
\end{theorem}

\begin{remark}
A minimal index of complementary is bounded for log Del Pezzo surfaces with
standard coefficients \cite[7.1]{Sh1}.
A hypothesis is that this index is not
more then 66.
This implies that we can take $\QQ$-complement of special type for four
dimensional singularities. Namely,
let $d_i=r_i/q_i$ then $q_i\le 66$. It allows to consider the finite number
variants of $D$.
\end{remark}

\begin{conjecture}
Let $(X\ni P)$ be an $n$-dimensional $(n\ge 3)$ canonical (log canonical)
hypersurface singularity. Let $f$ be a log resolution and
$\min\{a(E,0)| f(E)=P\}\ge n-2$. Then there exists a hyperplane section
$H$ such that $(X,H)$ is plt (lc). Hence $(X\ni P)$ is not weakly exceptional
(not exceptional) singularity \cite[4.8]{PrPlt}.
\end{conjecture}

\end{document}